\documentclass{gtart_h}  

\def\ifplaintex{\expandafter\ifx\csname documentclass\endcsname\relax}

\def\ifplaintex{\expandafter\ifx\csname documentclass\endcsname\relax}

\ifplaintex 
\else
\fi

\expandafter\ifx\csname epsfbox\endcsname\relax\input epsf\fi

\def\gt{{\mathsurround=0pt\it $\cal G\mskip-2mu$eometry \&\ 
$\cal T\!\!$opology}}        

\def\gtp{{\mathsurround=0pt\it $\cal G\mskip-2mu$eometry \&\ 
$\cal T\!\!$opology $\cal P\!$ublications}}  

\def\lognumber#1{\def\thelognumber{#1}}
\def\volumenumber#1{\def\thevolumenumber{#1}}
\def\papernumber#1{\def\thepapernumber{#1}}
\def\volumeyear#1{\def\thevolumeyear{#1}}

\def\pagenumbers#1#2{\def\startpage{#1}\def\finishpage{#2}}
\def\published#1{\def\publishdate{#1}}
\def\proposed#1{\def\theproposer{#1}}
\def\seconded#1{\def\theseconders{#1}}
\def\received#1{\def\receiveddate{#1}}

\def\accepted#1{\def\accepteddate{#1}}
\def\asciititle#1{\def\theasciititle{#1}}
\def\covertitle#1{\def\thecovertitle{#1}}
\def\coverauthors#1{\def\thecoverauthors{#1}}
\def\asciiauthors#1{\def\theasciiauthors{#1}}
\def\asciiaddress#1{\def\theasciiaddress{#1}}
\def\asciiemail#1{\def\theasciiemail{#1}}

\long\def\asciiabstract#1{\long\def\theasciiabstract{#1}}
\def\asciikeywords#1{\def\theasciikeywords{#1}}

\let\\\par\let\thelognumber\relax
\let\thevolumenumber\relax\let\thepapernumber\relax
\let\thevolumeyear\relax\let\thesamplenumber\relax\let\startpage\relax
\let\finishpage\relax\let\publishdate\relax\let\receiveddate\relax
\let\reviseddate\relax\let\accepteddate\relax\let\theasciititle\relax
\let\thecovertitle\relax\let\theasciiauthors\relax\let\theasciiaddress\relax
\let\theasciiabstract\relax\let\theasciikeywords\relax
\let\theasciiemail\relax\let\theshortauthors\relax\let\theshorttitle\relax
\let\thecoverauthors\relax

\long\def\maketitlep{   

\count0=\startpage

\gt\hfill      
\hbox to 77pt{\vbox to 0pt{\vglue -15pt\epsfbox{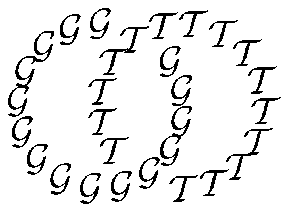}\vss}\hss}
\break
{\small\ifx\thesamplenumber\relax 
Volume \else Sample
\fi\thevolumenumber\ (\thevolumeyear)
\startpage--\finishpage\nl
Published: \publishdate}
\vglue 0.5truein plus 0.4fil minus 0.1truein

{\parskip=0pt\leftskip 0pt plus 1fil\def\\{\par\smallskip}{\ifplaintex\large
\else\Large\fi\bf\thetitle}\par\medskip}   

\vglue 0pt plus 0.1fil 

{\parskip=0pt\leftskip 0pt plus 1fil\def\\{\par}{\sc\theauthors}
\par\medskip}

\vglue 0pt plus 0.1fil 

{\small\parskip=0pt\let\newline\\
{\leftskip 0pt plus 1fil\def\\{\par}{\sl\theaddress}\par}
\expandafter\ifx\theemail\relax    
\relax\else\vglue 5pt plus 0.02fil minus 2pt\def\\{\stdspace{\rm 
and}\stdspace} 
\cl{Email:\stdspace\tt\theemail}\fi
\ifx\theurl\relax                  
\relax\else\vglue 5pt plus 0.02fil minus 2pt\def\\{\stdspace{\rm 
and}\stdspace}
\cl{URL:\stdspace\tt\theurl}\fi\par}

\vglue 7pt plus 0.3fil minus 3pt

{\bf Abstract}
\vglue 5pt plus 0.1fil minus 2pt

\theabstract

\vglue 7pt plus 0.3fil minus 3pt

{\bf AMS Classification numbers}\quad Primary:\quad \theprimaryclass

Secondary:\quad \thesecondaryclass

\vglue 5pt plus 0.3fil minus 2pt

{\bf Keywords:}\quad \thekeywords

\vglue 10pt plus 0.5fil minus 5pt

{\small  Proposed: \theproposer\hfill Received: \receiveddate\nl
Seconded: \theseconders\hfill 
\ifx\reviseddate\relax                         
Accepted: \accepteddate                        
\else
Revised: \reviseddate                          
\fi}
\eject
}       

\font\phead=cmsl9 scaled 950
\font=cmsl9 scaled 1050
\font\pnum=cmbx10 scaled 913
\font\lnum=cmbx10 
\font\pfoot=cmsl9 scaled 950
\font=cmsl9 scaled 1050
\ifplaintex
\headline{\vbox to 0pt{\vskip -4.5mm\line{\small\phead\ifnum
\count0=\startpage ISSN 1364-0380 (on line)
1465-3060 (printed) \hfill {\pnum\folio}\else\ifodd\count0\def\\{ }%
\ifx\theshorttitle\relax\thetitle\else\theshorttitle\fi\hfill{\pnum\folio}
\else\def\\{ and }{\pnum\folio}\hfill\ifx\theshortauthors\relax\theauthors
\else\theshortauthors\fi\fi\fi}\vss}}
\footline{\vbox to 0pt{\vglue 0mm\line{\small\pfoot\ifnum\count0=\startpage
\copyright\ \gtp\hfill\else
\gt, Volume \thevolumenumber\ (\thevolumeyear)\hfill\fi}\vss
}}
\else
\makeatletter
\def\@oddhead{{\small\lhead\ifnum\count0=\startpage ISSN 1364-0380 (on line)
1465-3060 (printed) \hfill {\lnum\number\count0}\else\ifodd\count0
\def\\{ }\ifx\theshorttitle\relax \thetitle \else\theshorttitle\fi\hfill
{\lnum\number\count0}\else\def\\{ and }{\lnum\number\count0}
\hfill\ifx\theshortauthors\relax 
\theauthors\else\theshortauthors\fi\fi\fi}}\def\@evenhead{\@oddhead}
\def\@oddfoot{\small\lfoot\ifnum\count0=\startpage\copyright\ \gtp\hfill\else
\gt, Volume \thevolumenumber\ (\thevolumeyear)\hfill\fi}
\def\@evenfoot{\@oddfoot}
\makeatother
\fi

\newwrite\gtoutfile
\long\gdef\makeheadfile{  
{\def\\{, }\def\s{ }
\immediate\openout\gtoutfile head.xxx
\immediate\write\gtoutfile{Proxy-for: \ifx\theasciiauthors\relax
\theauthors\else\theasciiauthors\fi\s<\ifx\theasciiemail\relax\theemail\else\theasciiemail\fi>}
\immediate\write\gtoutfile{\noexpand\\}
\immediate\write\gtoutfile{Authors: \ifx\theasciiauthors\relax
\theauthors\else\theasciiauthors\fi}
{\def\\{ }\immediate\write\gtoutfile{Title: \ifx\theasciititle\relax
\thetitle\else\theasciititle\fi}}
\immediate\write\gtoutfile{Subj-class: GT or SG or MG etc}
\immediate\write\gtoutfile{MSC-class: \theprimaryclass\ifx\thesecondaryclass\relax\else, \thesecondaryclass\fi}
\immediate\write\gtoutfile{Journal-ref: Geom. Topol. \thevolumenumber
(\thevolumeyear) \startpage-\finishpage}
\immediate\write\gtoutfile{Comments: Published by Geometry and Topology at}
\immediate\write\gtoutfile{\s\s http://www.maths.warwick.ac.uk/gt/GTVol\thevolumenumber/paper\thepapernumber.abs.html}
\immediate\write\gtoutfile{\noexpand\\}
\immediate\write\gtoutfile{}
\ifx\theasciiabstract\relax
\immediate\write\gtoutfile{\theabstract}\else
\immediate\write\gtoutfile{\theasciiabstract}\fi
\immediate\write\gtoutfile{}
\immediate\write\gtoutfile{\noexpand\\}
\immediate\write\gtoutfile{}
\immediate\closeout\gtoutfile}}  

\def\maketitlepage{\maketitlep\makeheadfile}
\let\maketitle\maketitlepage

\lognumber{426}
\volumenumber{8}\papernumber{24}\volumeyear{2004}
\pagenumbers{925}{945}
\received{21 February 2004}
\published{9 June 2004}
\accepted{29 May 2004}
\proposed{Peter Ozsv\'ath}
\seconded{John Morgan, Tomasz Mrowka}

\usepackage{amsmath, amssymb, graphicx, graphs, verbatim}

\newtheorem{thm}{Theorem}[section]
\newtheorem{cor}[thm]{Corollary}

\newtheorem{prop}[thm]{Proposition}

\theoremstyle{definition}

\theoremstyle{remark}

\numberwithin{equation}{section}

\newcommand{\al}{\alpha}

\newcommand{\ga}{\gamma}

\newcommand{\La}{\Lambda}

\newcommand{\om}{\omega}

\newcommand{\Si}{\Sigma}

\newcommand{\Dm}{{\mathbb {D}}_m}

\newcommand{\x}{\times}
\newcommand{\s}{\mathbf s}

\renewcommand{\t}{\mathbf t}

\newcommand{\C}{\mathbb C}
\newcommand{\Z}{\mathbb Z}
\newcommand{\N}{\mathbb N}
\newcommand{\Q}{\mathbb Q}

\newcommand{\del}{\partial}

\newcommand{\lra}{\longrightarrow}
\newcommand{\hf}{{{\widehat {HF}}}}

\newcommand{\Li}{\mathbb {L}}

\newcommand{\mk}{{\overline {K}}}

\DeclareMathOperator{\Hom}{Hom}
\DeclareMathOperator{\tb}{tb}
\DeclareMathOperator{\TB}{TB}

\DeclareMathOperator{\PD}{PD}
\DeclareMathOperator{\Spin}{Spin}

\DeclareMathOperator{\Ima}{Im}
\begin{document}

\title{Ozsv\'ath--Szab\'o invariants and tight\\contact three--manifolds, I}
\covertitle{Ozsv\noexpand\'ath--Szab\noexpand\'o invariants and tight\\contact three--manifolds, I}
\asciititle{Ozsvath-Szabo invariants and tight contact three-manifolds, I}

\author{Paolo Lisca\\Andr\'{a}s I Stipsicz}
\coverauthors{Paolo Lisca\\Andr\noexpand\'{a}s I Stipsicz}
\asciiauthors{Paolo Lisca\\Andras I Stipsicz}

\address{Dipartimento di Matematica, Universit\`a 
di Pisa \\I-56127 Pisa, ITALY} 
\secondaddress{R\'enyi Institute of Mathematics, Hungarian 
Academy of Sciences\\H-1053 Budapest, Re\'altanoda utca 13--15, Hungary}

\asciiaddress{Dipartimento di Matematica, Universita 
di Pisa \\I-56127 Pisa, ITALY\\and\\Renyi Institute of Mathematics, Hungarian 
Academy of Sciences\\H-1053 Budapest, Realtanoda utca 13--15, Hungary}

\gtemail{\mailto{lisca@dm.unipi.it}{\rm \qua 
and\qua}\mailto{stipsicz@math-inst.hu}}

\asciiemail{lisca@dm.unipi.it, stipsicz@math-inst.hu}

\begin{abstract}
Let $S^3_r(K)$ be the oriented 3--manifold obtained by rational
$r$--surgery on a knot $K\subset S^3$. Using the contact
Ozsv\'ath--Szab\'o invariants we prove, for a class of knots $K$
containing all the algebraic knots, that $S^3_r(K)$ carries positive,
tight contact structures for every $r\neq 2g_s(K)-1$, where $g_s(K)$
is the slice genus of $K$. This implies, in particular, that the
Brieskorn spheres $-\Sigma (2,3,4)$ and $-\Sigma (2,3,3)$ carry tight,
positive contact structures. As an application of our main result we
show that for each $m\in\N$ there exists a Seifert fibered rational
homology 3--sphere $M_m$ carrying at least $m$ pairwise
non--isomorphic tight, nonfillable contact structures.
\end{abstract}
\asciiabstract{%
Let S^3_r(K) be the oriented 3--manifold obtained by rational
r-surgery on a knot K in S^3. Using the contact Ozsvath-Szabo
invariants we prove, for a class of knots K containing all the
algebraic knots, that S^3_r(K) carries positive, tight contact
structures for every r not= 2g_s(K)-1, where g_s(K) is the slice genus
of K. This implies, in particular, that the Brieskorn spheres
-Sigma(2,3,4) and -Sigma(2,3,3) carry tight, positive contact
structures. As an application of our main result we show that for each
m in N there exists a Seifert fibered rational homology 3-sphere M_m
carrying at least m pairwise non-isomorphic tight, nonfillable contact
structures.}

\primaryclass{57R17} 
\secondaryclass{57R57} 
\keywords{Tight, fillable
contact structures, Ozsv\'ath--Szab\'o invariants}
\asciikeywords{Tight, fillable
contact structures, Ozsvath-Szabo invariants}

\maketitlepage

\section{Introduction}\label{s:intro}
According to a classical result of Lutz and Martinet, every closed,
oriented 3--manifold admits a positive contact structure. In fact,
every oriented 2--plane field on an oriented 3--manifold is homotopic
to a positive contact structure. The proof of the Lutz--Martinet
theorem --- relying on contact surgery along transverse links in the
standard contact 3-sphere~\cite{Geigb} --- typically produces
overtwisted contact structures.  (For a proof of the Lutz--Martinet
theorem using contact surgery along Legendrian links see \cite{DGS}.)
Finding~\emph{tight} contact structures on a closed 3--manifold is, in
general, much more difficult, indeed impossible for the Poincar\'e
homology 3--sphere with its natural orientation reversed~\cite{EH2}.

Let $Y$ be a closed, oriented 3--manifold. Consider the following
problem:
\begin{enumerate}
\item[(P)]
Does $Y$ carry a positive, tight contact structure?
\end{enumerate}
Until recently, the two most important methods to deal with problem
(P) were Eliashberg's~\emph{Legendrian surgery} as used eg~by Gompf
in~\cite{G}, and the~\emph{state traversal} method, developed by Ko
Honda and based on Giroux's theory of convex surfaces. The limitations
of these two methods come from the fact that Legendrian surgery can
only prove tightness of Stein fillable contact structures, while the
state traversal becomes combinatorially unwieldy in the absence of
suitable incompressible surfaces. For example, both methods fail to
deal with problem (P) when $Y$ is one of the Brieskorn spheres
$-\Sigma (2,3,4)$ or $-\Sigma (2,3,3)$, because these Seifert fibered
3--manifolds do not contain vertical incompressible tori, nor do they
carry symplectically fillable contact structures~\cite{PLpos}.

The purpose of the present paper is to show that 
contact Ozsv\'ath--Szab\'o invariants~\cite{OSz6} can be effectively
combined with contact surgery~\cite{DG1, DG2} to tackle problem
(P). In particular, it follows from Theorem~\ref{t:main} below that
$-\Sigma (2,3,4)$ and $-\Sigma (2,3,3)$ do indeed carry tight,
positive contact structures. Moreover, such contact structures admit an
explicit description (cf~Corollary~\ref{c:brieskorn} and the
following remark).

In order to state our main result we need to introduce some notation.
Recall that the~\emph{standard contact structure} on $S^3$ is the
2--dimensional distribution $\xi_{\rm st}\subset TS^3$ given by the
complex tangents, where $S^3$ is viewed as the boundary of the unit
4--ball in $\C^2$. We say that a knot in $S^3$ is~\emph{Legendrian} if
it is everywhere tangent to $\xi_{\rm st}$. To every Legendrian knot
$L\subset S^3$ one can associate its Thurston--Bennequin
number~$\tb(L)\in\Z$, which is invariant under Legendrian isotopies of
$L$~\cite{Be}. Given a knot $K\subset S^3$, let $\TB(K)$ denote
the~\emph{maximal} Thurston--Bennequin number of $K$, defined as
\begin{enumerate}
\item[]
$\TB(K)=\max \{ tb(L) \mid L$ is Legendrian and smoothly isotopic to $K\}$.
\end{enumerate}
Let $g_s(K)$ denote the~\emph{slice genus} (aka~the~\emph{4--ball
genus}) of $K$. Let $S^3_r(K)$ be the oriented 3--manifold given by
rational $r$--surgery on a knot $K\subset S^3$.

\begin{thm}\label{t:main}
Let $K\subset S^3$ be a knot such that 
\[
g_s(K)>0\quad\text{and}\quad\TB(K)=2g_s(K)-1. 
\]
Then, the oriented 3--manifold $S^3_r(K)$ carries positive, tight contact
structures for every $r\neq 2g_s(K)-1$. 
\end{thm}

\rk{Remark} By the slice Bennequin inequality~\cite{Ru}, for any knot
$K\subset S^3$ we have
\[
\TB(K)\leq 2g_s(K)-1.
\]
Moreover, by~\cite{BW, Br} (see~\cite[page~123]{Be}), if $K$ is an
algebraic knot then
\[
\TB(K)=2g(K)-1,
\]
where $g(K)$ is the Seifert genus of $K$. Since $g_s(K)\leq g(K)$, it
follows that the family of knots $K$ satisfying the assumption of
Theorem~\ref{t:main} contains all nontrivial algebraic knots. In fact,
there are non--fibered, hence non--algebraic, knots satisfying the
same assumption, as for example certain negative twist knots.%
\footnote{Let $K_q$ be a twist knot with $q<0$ twists
(cf~\cite[page~112]{Ro}). It is easy to find a Legendrian
representative of $K_q$ with Thurston--Bennequin number equal to
$1$. On the other hand, by resolving the clasp it follows that
$g_s(K_q) \leq 1$.  Therefore the slice Bennequin inequality implies
$g_s(K_q)=TB(K_q)=1$. The knots $K_q$ are not fibered for $q<-1$
because the leading coefficient of their Alexander polynomial is not
equal to $1$.}

Let $T\subset S^3$ be the right--handed trefoil. Since $T$ is
algebraic, Theorem~\ref{t:main} applies. In particular, since
$S^3_2(T)=-\Sigma (2,3,4)$ and $S^3_3(T)=-\Sigma (2,3,3)$,
Theorem~\ref{t:main} immediately implies the following result, which
solves a well--known open problem~\cite[Question 8]{Et}:

\begin{cor}\label{c:brieskorn}
The Brieskorn spheres $-\Sigma (2,3,3)$ and $-\Sigma (2,3,4)$ carry
positive, tight contact structures. \qed
\end{cor}

\rk{Remarks} (1)\qua The proof of Theorem~\ref{t:main} shows that
Figures~\ref{f:knot} and~\ref{f:link} below provide explicit
descriptions of the tight contact structures of
Corollary~\ref{c:brieskorn}.

(2)\qua Theorem~\ref{t:main} is optimal for the right--handed trefoil knot
 $T=T_{3,2}$, because $S^3_1(T)=-\Si(2,3,5)$ is known not to carry
 positive, tight contact structures~\cite{EH2}. On the other hand, it
 is natural to ask whether the same is true for other torus knots. We
 address this question in the companion paper~\cite{LS4}.

Recall that a~\emph{symplectic filling} of a contact three--manifold
$(Y,\xi)$ is a pair $(X,\om)$ consisting of a smooth, compact,
connected four--manifold $X$ and a symplectic form $\om$ on $X$ such
that, if $X$ is oriented by $\om\wedge\om$, $\del X$ is given the
boundary orientation and $Y$ is oriented by $\xi$, then $\del X=Y$ and
$\om|_{\xi}\neq 0$ at every point of $\del X$.  As an application of
Theorem~\ref{t:main} we prove the following result, which should be
compared with the results of~\cite{LS2, LS3}.

\begin{thm}\label{t:nonfill}
For each $m\in \N$ there is a Seifert fibered rational homology sphere
$M_m$ carrying at least $m$ pairwise non--isomorphic tight,
not symplectically fillable contact structures.
\end{thm}

The paper is organized as follows. In Section~\ref{s:surgeries} we
describe the necessary background in contact surgery and Heegaard
Floer theory.  In Sections~\ref{s:main} and~\ref{s:nonfill} we prove,
respectively, Theorems~\ref{t:main} and~\ref{t:nonfill}.

\rk{Acknowledgements}The first author was partially supported by
  MURST, and he is a member of EDGE, Research Training Network
  HPRN-CT-2000-00101, supported by The European Human Potential
  Programme. The second author would like to thank Peter Ozsv\'ath and
  Zolt\'an Szab\'o for many useful discussions regarding their joint
  work. The second author was partially supported by OTKA T34885.

\section{Surgeries and Ozsv\'ath--Szab\'o invariants}\label{s:surgeries}

\sh{Contact surgery}

Let $(Y, \xi)$ be a contact 3--manifold. The framing of a
Legendrian knot $K\subset Y$ naturally induced by $\xi$ is called
the~\emph{contact framing} of $K$. Given a Legendrian knot $K$ in a
contact 3--manifold $(Y,\xi)$ and a non--zero rational number
$r\in\Q$, one can perform \emph{contact $r$--surgery} along $K$ to
obtain a new contact 3--manifold $(Y',\xi')$~\cite{DG1, DG2}. Here
$Y'$ is the 3--manifold obtained by smooth $r$--surgery along $K$
with respect to the contact framing, while $\xi'$ is constructed by
extending $\xi$ from the complement of a standard neighborhood of $K$
to a tight contact structure on the glued--up solid torus.  If $r\neq
0$ such an extension always exists, and for $r=\frac 1k$ $(k\in \Z)$
it is unique \cite{H2}.  When $r=-1$ the corresponding contact surgery
coincides with Legendrian surgery along $K$~\cite{el,G,W}.

\begin{figure}[ht]
\begin{center}\small
\setlength{\unitlength}{1mm}
\begin{picture}(41,30)
\put(33,25){$+1$}
\includegraphics[height=3cm]{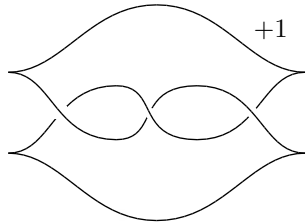}
\end{picture}    
\end{center}
\caption{\quad A contact structure on $-\Sigma (2,3,4)$}
\label{f:knot}
\end{figure}

As an illustration of the contact surgery construction, consider the
Legendrian trefoil knot $\mathcal T$ represented by the Legendrian
front (see eg~\cite{G} for notation) of Figure~\ref{f:knot}. Since
the coefficient $+1$ represents the contact surgery coefficient and
$\tb({\mathcal T})=1$, the picture represents a contact structure on
the oriented 3--manifold obtained by a smooth $(+2)$--surgery on a
right--handed trefoil knot, that is on $-\Si(2,3,4)$.

According to~\cite[Proposition~7]{DG2}, a contact
$r=\frac{p}{q}$--surgery ($p,q\in \N$) on a Legendrian knot $K$ is
equivalent to a contact $\frac{1}{k}$--surgery on $K$ followed by a
contact $\frac{p}{q-kp}$--surgery on a Legendrian pushoff of $K$ for
any integer $k\in \N$ such that $q-kp< 0$. Moreover,
by~\cite[Proposition~3]{DG2} any contact $r$--surgery along $K\subset
(Y, \xi )$ with $r<0$ is equivalent to a Legendrian surgery along a
Legendrian link $\Li =\cup_{i=0}^m L_i$ which is determined via a
simple algorithm by the Legendrian knot $K$ and the contact surgery
coefficient $r$. The algorithm to obtain $\Li $ is the following. Let
\[
[a_0, \ldots ,a_m],\quad a_0,\ldots a_m\geq 2 
\]
be the continued fraction expansion of $1-r$. To obtain the first
component $L_0$, push off $K$ using the contact framing and stabilize
it $a_0-2$ times. Then, push off $L_0$ and stabilize it $a_1-2$
times.  Repeat the above scheme for each of the remaining pivots of
the continued fraction expansion.  Since there are $a_i-1$
inequivalent ways to stabilize a Legendrian knot $a_i-2$ times, this
construction yields $\Pi_{i=0}^m (a_i-1)$ potentially different
contact structures.

\begin{figure}[ht]
\setlength{\unitlength}{1mm}
\begin{center}\small
\begin{picture}(118,38)
\put(0,0){\includegraphics[height=3.8cm]{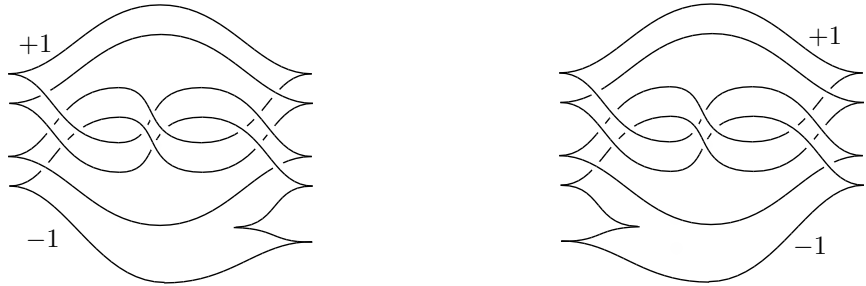}}
\put(3,31){$+1$}
\put(108,32){$+1$}
\put(4,5){$-1$}
\put(106,4){$-1$}
\end{picture}    
\end{center}
\caption{\quad Contact structures on $-\Sigma (2,3,3)$}
\label{f:link}
\end{figure}

For example, according to the algorithm just described, any contact
$(+2)$--surgery on $\mathcal T$ is equivalent to one of
the contact surgeries of Figure~\ref{f:link} (the coefficients
indicate surgery with respect to the contact framings).

Since, by~\cite[Proposition~9]{DG1}, a contact $\frac{1}{k}$--surgery
($k\in \N$) on a Legendrian knot $K$ can be replaced by $k$ contact
$(+1)$--surgeries on $k$ Legendrian pushoffs of $K$, it follows that
any contact rational $r$--surgery ($r\neq 0$) can be replaced by
contact $(\pm 1)$--surgery along a Legendrian link; for a related
discussion see also~\cite{DGS, LS3}.

\sh{The Ozsv\'ath--Szab\'o invariants of 3--manifolds}

The Ozsv\'ath--Szab\'o invariants~\cite{OSzF1, OSzF2, OSzF4} assign to
each oriented spin$^c$ 3--manifold $(Y, \s)$ a finitely generated
Abelian group $\hf (Y, \s)$, and to each oriented spin$^c$ cobordism
$(W, \t)$ between $(Y_1,\s_1)$ and $(Y_2,\s_2)$ a homomorphism
\[
F_{W, \t}\colon \hf (Y_1, \s _1) \to \hf (Y_2, \s _2).
\]

For simplicity, in the following we will use these homology theories with
$\Z/2\Z$ coefficients. In this setting, $\hf(Y, \s )$ is a finite
dimensional vector space over the field $\Z/2\Z$. Define
\[
\hf(Y) = \bigoplus_{\s\in \Spin ^c (Y)} \hf (Y , \s).
\] 

Since there are only finitely many spin$^c$ structures with
nonvanishing invariants \cite[Theorem~7.1]{OSzF2}, $\hf (Y)$ is still
finite dimensional.

An important ingredient of our proofs is the following result, which
appears implicitly in the papers of Ozsv\'ath and Szab\'o (see
especially~\cite{OSz4genus}). We provide a detailed proof for completeness.

\begin{prop}\label{p:adjunction}
Let $W$ be a cobordism obtained by attaching a 2--handle to a
3--manifold $Y$ with $b_1(Y)=0$. Let $\t_0\in\Spin^c(W)$, and
suppose that $W$ contains a smoothly embedded, closed, oriented
surface $\Sigma$ of genus $g(\Si)>0$ such that
\[
\Si\cdot\Si\geq 0\quad\text{and}\quad
|\langle c_1(\t_0),[\Si]\rangle| + \Si\cdot\Si > 2g(\Si)-2. 
\]
Then, $F_{W,\t_0}=0$.
\end{prop}

\begin{proof}
Arguing by contradiction, suppose that $F_{W,\t_0}\neq 0$. By a
fundamental property of the invariants~\cite{OSzF4} there are only
finitely many spin$^c$ structures $\t_1,\ldots,\t_k\in\Spin^c(W)$
such that $F_{W,{\t_i}}\neq 0$. Moreover, by~\cite[Theorem~3.6]{OSzF4}
we have
\[
F_{W,\t_0}\neq 0\quad\Longleftrightarrow\quad F_{W,{\overline\t_0}}\neq 0
\]
where ${\overline\t_0}$ is the spin$^c$ structure conjugate to
$\t_0$. Therefore, up to replacing $\t_0$ with one of the $\t_i$'s we
may assume that 
\begin{equation}\label{e:max}
\langle c_1(\t_0), [\Si]\rangle=
|\langle c_1(\t_0), [\Si]\rangle|=
\max\{\langle c_1(\t_i), [\Si]\rangle\ |\ i=1,\ldots,k\}.
\end{equation}
Let $\Si\cdot\Si=n$, and let $\widehat W$ be the smooth 4--manifold
obtained by blowing up $W$ at $n$ distinct points of
$W\setminus\Si$. Choose exceptional classes
\[
e_1,\ldots,e_n\in H_2(\widehat W)
\]
and let $\widehat\t_0$ denote the unique spin$^c$ structure on $\widehat
W$ such that $\widehat\t_0|_W=\t$ and $\langle c_1(\widehat\t_0),
e_i\rangle =1$ for every $i=1,\ldots, n$.

Let $\widehat\Si\subset\widehat W$ be a smooth, oriented surface 
obtained by piping $\Si$ to the $n$ exceptional spheres, so that 
\[
[\widehat\Si] = [\Si] + \sum_{i=1}^n e_i.
\]
Let $\ga\subset\widehat W$ be a properly embedded arc (disjoint from $Y$
and $\widehat\Si$ away from its endpoints) connecting $Y$ to
$\widehat\Si$. Denote by $\widehat W_1$ a closed regular neighborhood
of the union $Y\cup\ga\cup\widehat\Si$, and let $\widehat W_2$ be the
closure of $\widehat W\setminus\widehat W_1$.

Let
\[
{\mathcal S}=\left\{\widehat\t\in\Spin^c(\widehat W)\qua|\qua\widehat\t
\vert_{\widehat W_i}=\widehat\t_0\vert_{\widehat W_i},\ i=1,2\right\}. 
\]
By the composition law~\cite[Theorem~3.4]{OSzF4} we have
\begin{equation}\label{e:composition}
F_{\widehat W_2,\widehat\t_0\vert_{\widehat W_2}}\circ 
F_{\widehat W_1,\widehat\t_0\vert_{\widehat W_1}}=
\sum_{\widehat\t\in\mathcal S} F_{\widehat W,\widehat\t}.
\end{equation}
We are going to show that the sum at the right hand side
of~\eqref{e:composition} admits at most one nontrivial term. 
In fact, we shall prove that  
\[
\widehat\t\in\mathcal S\quad\text{and}\quad
F_{\widehat W,\widehat\t} \neq 0\quad\Longrightarrow\quad
\widehat\t=\widehat\t_0. 
\]
Recall that $Spin^c({\widehat {W}})$ admits a free and transitive
action of $H^2({\widehat {W}}; \Z )$. Hence, there is an element $L\in
H^2 ({\widehat {W}}; \Z )$ such that
\[
\widehat\t - \widehat\t_0 = L.
\]
Since 
\[
\widehat\t \vert_{\widehat W_i}=\widehat\t_0\vert_{\widehat
W_i},\quad i=1,2, 
\]
we have, in particular, $L\vert_Y=0$. Therefore $L$ is the image of an
element $A\in H^2 ({\widehat W},Y;\Z)$ under the restriction map $H^2
({\widehat {W}},Y;\Z)\to H^2 ({\widehat W};\Z)$. Our plan is to show
that $\widehat\t=\widehat\t_0$ by proving that $A=0$. Since
\[
H_1(W,Y;\Z )\cong H_1 ({\widehat {W}},Y;\Z )=0, 
\]
the universal coefficient theorem implies that 
\[
H^2({\widehat W}, Y;\Z )\cong 
\Hom (H_2({\widehat W}, Y; \Z ), \Z ), 
\]
therefore to show $A=0$ it is enough to show $2A=0$, and $2A$ is
determined by its values on the elements of $H_2({\widehat W}, Y;\Z
)$. But since $b_1(Y)=0$, it suffices to show that $2A$ evaluates
trivially on the image of the map 
\[
i_*\co H_2({\widehat W};\Z)\lra H_2({\widehat W}, Y;\Z).
\]
On the other hand, since $\widehat\Si\subset\widehat W_1$, if
$\widehat\t\in{\mathcal S}$ then $\langle c_1(\widehat\t),
[\widehat\Si] \rangle = \langle c_1(\widehat\t_0), [\widehat\Si]
\rangle$, ie,
\begin{equation}\label{e:evaluation}
\langle c_1(\widehat\t|_W),[\Si]\rangle + 
\sum_{i=1}^n \langle c_1(\widehat\t), e_i\rangle = 
\langle c_1(\t_0),[\Si]\rangle + n.
\end{equation}
Moreover, by the blow--up formula~\cite[Theorem~3.7]{OSzF4} if
$\widehat\t\in\Spin^c(\widehat W)$ then
\[
F_{W,\t|_W}\neq 0\quad\Longleftrightarrow\quad 
F_{{\widehat W},{\widehat\t}}\neq 0
\quad\Longrightarrow\quad |\langle c_1(\widehat\t), e_i\rangle|=1,
\quad i=1,\ldots, n.
\]
Therefore, if $F_{{\widehat W},{\widehat\t}}\neq 0$, by 
Equations~\eqref{e:max} and~\eqref{e:evaluation} we have
\[
\langle c_1(\widehat\t|_W),[\Si]\rangle = 
\langle c_1(\t_0),[\Si]\rangle\quad\text{and}\quad
\langle c_1(\widehat\t), e_i\rangle = 
\langle c_1(\widehat\t_0), e_i\rangle=1,\quad i=1,\ldots,n.
\]
It follows that $c_1(\widehat\t)=c_1(\widehat\t_0)$. Therefore, 
for every $\al\in H_2({\widehat W};\Z)$ we have
\[
\langle 2A, i_*(\al)\rangle = \langle 2L, \al\rangle = 
\langle c_1(\widehat\t) - c_1(\widehat\t_0), \al\rangle = 0.
\]
Thus, $\widehat\t=\widehat\t_0$, and the right--hand side of
Equation~\eqref{e:composition} reduces to $F_{\widehat
W,\widehat\t_0}$.

Now observe that $\widehat W_1$ is a cobordism from $Y$ to
$Y\# S^1\x\widehat\Si$, and since
\[
\langle c_1(\widehat\t_0), [\widehat\Si]\rangle =
\langle c_1(\t_0),[\Si]\rangle + n >
2g(\widehat\Si) - 2,
\]
by the adjunction inequality~\cite[Theorem~7.1]{OSzF2} the group
\[
\hf(Y\# S^1\x\widehat\Si, \widetilde\t_0\vert_{S^1\x\widehat\Si})
\]
is trivial. But this group is the domain of the map $F_{\widehat
W_2,\widehat\t_0\vert_{\widehat W_2}}$. Thus,
Equation~\eqref{e:composition} implies that $F_{\widehat
W,\widehat\t_0}=0$ and therefore $F_{W,\t _0}=0$, which gives the
desired contradiction.
\end{proof} 


\sh{Contact Ozsv\'ath--Szab\'o  invariants}

In~\cite{OSz6} Ozsv\'ath and Szab\'o defined an invariant 
\[
c(Y, \xi )\in \hf (-Y, \s _{\xi})/\langle \pm 1 \rangle
\] 
for a contact 3--manifold $(Y,\xi )$, where $\s _{\xi }$ denotes the
spin$^c$ structure induced by the contact structure $\xi$. Since in
this paper we are using this homology theory with $\Z/2\Z$
coefficients, the above sign ambiguity for $c(Y,\xi)$ does not occur.
It is proved in~\cite{OSz6} that if $(Y, \xi)$ is overtwisted then
$c(Y,\xi)=0$, and if $(Y,\xi)$ is Stein fillable then $c(Y,\xi)\neq
0$. In particular, $c(S^3,\xi_{\rm st})\neq 0$. We are going to use the
properties of $c(Y,\xi)$ described in the following theorem and
corollary.

\begin{thm}[\cite{LS3}, Theorem~2.3]\label{t:+1}
Suppose that $(Y',\xi ')$ is obtained from $(Y,\xi )$ by a contact
$(+1)$--surgery. Let $-X$ be the cobordism induced by the surgery with
reversed orientation. Define 
\[
F_{-X} := \sum_{\t\in\Spin^c(-X)} F_{-X,\t}.
\]
Then,
\[
F_{-X} (c(Y , \xi ))= c(Y',\xi '). 
\]
In particular, if $c(Y',\xi')\neq 0$ then $(Y,\xi)$ is tight. 
\qed\end{thm}

\begin{cor}[\cite{LS3}, Corollary~2.4]\label{c:surg}
If $c(Y_1,\xi _1)\neq 0$ and $(Y_2,\xi_2)$ is obtained from
$(Y_1,\xi_1)$ by Legendrian surgery along a Legendrian knot, then
$c(Y_2, \xi_2)\neq 0$. In particular, $(Y_2,\xi_2)$ is tight.
\qed\end{cor}

\sh{The surgery exact triangle}

Here we describe what is usually called the~\emph{surgery exact
triangle} for the Ozsv\'ath--Szab\'o homologies.

\begin{figure}[ht]
\setlength{\unitlength}{1mm}
\begin{center}\small
\begin{picture}(119,40)
\put(0,0){\includegraphics[height=4cm]{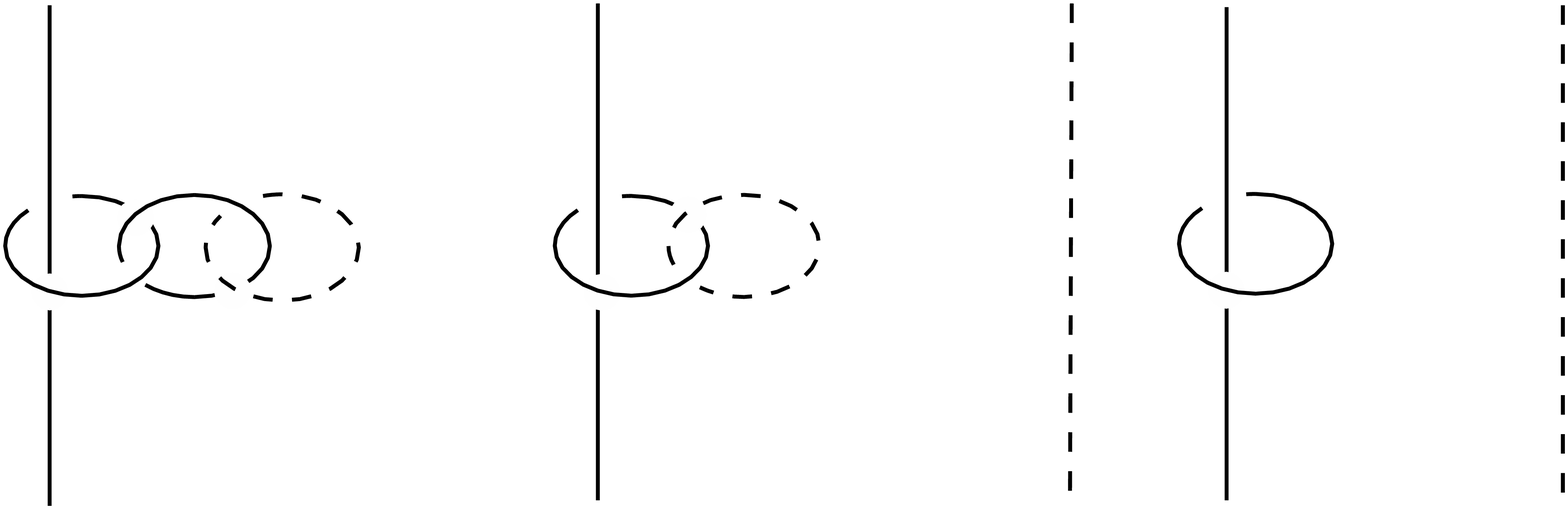}}
\put(6.5,4){$K$}
\put(7.5,14){$N$}
\put(15,14){$C$}
\put(23,14){$D$}
\put(45,4){$K$}
\put(79,4){$K$}
\put(90,4){$K$}
\put(114,4){$K$}
\put(7,35){$f$}
\put(45,35){$f$}
\put(79,35){$f$}
\put(90,35){$f$}
\put(114,35){$f$}
\put(32,20){$\cong$}
\put(67,20){$\cong$}
\put(102,20){$\cong$}
\put(-3,20){$-1$}
\put(12,26){$-1$}
\put(19,26){$-1$}
\put(46,26){$0$}
\put(54,26){$0$}
\put(92,26){$0$}
\end{picture}    
\end{center}
\caption{The diffeomorphism between $W$ and $X$}
\label{f:cobordism}
\end{figure}

Let $Y$ be a closed, oriented 3--manifold and let $K\subset Y$ be a
framed knot with framing $f$. Let $Y(K)$ denote the 3--manifold given
by surgery along $K\subset Y$ with respect to the framing $f$. The
surgery can be viewed at the 4--manifold level as a 4--dimensional
2--handle addition. The resulting cobordism $X$ induces a homomorphism
\[
F_X:=\sum_{\t\in\Spin^c(X)} F_{X,\t}\co \hf (Y)\to\hf(Y(K))
\]
obtained by summing over all spin$^c$ structures on $X$. Similarly,
there is a cobordism $U$ defined by adding a 2--handle to $Y(K)$ along
a normal circle $N$ to $K$ with framing $-1$ with respect to a normal
disk to $K$. The boundary components of $U$ are $Y(K)$ and the
3--manifold $Y'(K)$ obtained from $Y$ by a surgery along $K$ with
framing $f+1$. As before, $U$ induces a homomorphism
\[
F_U\co\hf(Y(K))\to\hf(Y'(K)).
\]
It is proved in~\cite[Theorem~9.16]{OSzF2}\footnote{In fact, the maps
$F_U$ and $F_X$ were defined in~\cite{OSzF2} by counting
pseudo--holomorphic triangles in a Heegaard triple, but an easy
comparison with the maps associated to 2--handles defined
in~\cite[Subsection~4.1]{OSzF4} shows that $F_U$ and $F_X$ are the
sums of maps associated to cobordisms given above (see the discussion
at the beginning of~\cite[Section~3]{abs}).} that
\begin{equation}\label{e:exact}
\ker F_U=\Ima F_X.
\end{equation}
The above construction can be repeated starting with $Y(K)$ and
$N\subset Y(K)$ equipped with the framing specified above: we get $U$
(playing the role previously played by $X$) and a new cobordism $V$
starting from $Y'(K)$, given by attaching a 4--dimensional 2--handle
along a normal circle $C$ to $N$ with framing $-1$ with respect to a
normal disk. It is easy to check that this last operation yields $Y$
at the 3--manifold level. Again, we have $\ker F_V=\Ima F_U$.
Moreover, we can apply the construction once again, and denote by $W$
the cobordism obtained by attaching a 2--handle along a normal circle
$D$ to $C$ with framing $-1$. In fact, $W$ is orientation--preserving
diffeomorphic to $X$. This fact is explained in
Figure~\ref{f:cobordism}, where the first picture represents $W$ and
the last picture represents $X$. In the figure, the framed dotted
circle is the attaching circle of the 2--handle. The first
diffeomorphism in Figure~\ref{f:cobordism} is obtained by ``blowing
down'' the framed knot $C$. In other words, the first two pictures
represent 2--handles attached to diffeomorphic 3--manifolds, and show
that the corresponding attaching maps commute with the given
diffeomorphism. The second diffeomorphism is obtained by a handle
slide, and the third diffeomorphism by erasing a cancelling pair.
It follows immediately from Equation~\eqref{e:exact} that the
homomorphisms $F_{X}, F_U$ and $F_V$ fit into the~\emph{surgery exact
triangle}:
\begin{equation}\label{e:triangle}
\begin{graph}(6,2)
\graphlinecolour{1}\grapharrowtype{2}
\textnode {A}(1,1.5){$\hf (Y)$}
\textnode {B}(5, 1.5){$\hf (Y(K))$}
\textnode {C}(3, 0){$\hf (Y'(K))$}
\diredge {A}{B}[\graphlinecolour{0}]
\diredge {B}{C}[\graphlinecolour{0}]
\diredge {C}{A}[\graphlinecolour{0}]
\freetext (3,1.8){$F_{X}$}
\freetext (4.6,0.6){$F_U$}
\freetext (1.4,0.6){$F_V$}
\end{graph}
\end{equation}

\rk{Remark}
Given an exact triangle of vector spaces and homomorphisms
\begin{equation*}
\begin{graph}(6,2)
\graphlinecolour{1}\grapharrowtype{2}
\textnode {A}(1,1.5){$V_1$}
\textnode {B}(5, 1.5){$V_2$}
\textnode {C}(3, 0){$V_3$}
\diredge {A}{B}[\graphlinecolour{0}]
\diredge {B}{C}[\graphlinecolour{0}]
\diredge {C}{A}[\graphlinecolour{0}]
\freetext (3,1.8){$F_3$}
\freetext (4.4,0.6){$F_1$}
\freetext (1.7,0.5){$F_2$}
\end{graph}
\end{equation*}
we have
\begin{equation}\label{e:triangleinequality}
\dim V_i \leq \dim V_j + \dim V_k
\end{equation}
for $\{i,j,k\}=\{1,2,3\}$. Moreover, equality holds
in~\eqref{e:triangleinequality} if and only if $F_i=0$.

\section{The proof of Theorem~\ref{t:main}}
\label{s:main}

Let $L$ be a Legendrian knot smoothly isotopic to $K$ with
\[
t:= \tb(L)=2g_s(K)-1. 
\]
Let $r\in\Q\setminus\{t\}$ and $r'=r-t$. Then, any contact
$r'$--surgery along $L$ yields a contact structure on $S^3_r(K)$.

If $r<t=2g_s(K)-1$ then $r'<0$. Since any contact $r'$--surgery along
$L$ can be realized by Legendrian surgery, the resulting contact
structure is Stein fillable and hence tight~\cite{Eli1}. Therefore, to
prove Theorem~\ref{t:main} it suffices to show that any contact $r'$--surgery
along $L$ with $r'>0$ yields a contact structure on $S^3_r(K)$ with
non--zero contact Ozsv\'ath--Szab\'o invariant.

Let $(Y_k,\xi_k)$, with $k$ any positive integer, denote the result of
contact $\frac 1k$--surgery along $L$. If $r'>0$, any contact
$r'$--surgery along $L$ is equivalent to a sequence of Legendrian
surgeries on $(Y_k,\xi_k)$ for some $k>0$. Therefore, by
Corollary~\ref{c:surg} it suffices to prove that the contact
invariants of $(Y_k,\xi _k)$ do not vanish. We claim that, for every
$k\geq 1$,
\begin{equation}\label{e:invariant}
\quad c(Y_k,\xi_k)\neq 0.
\end{equation}
We are going to prove the claim by induction on $k$. To start the
induction, we examine the case $k=1$ first.

Observe that $Y_1(L)=S^3_{2g_s}(K)$, and let $-X$ be the cobordism
induced by contact $(+1)$--surgery along $L$ with reversed
orientation. Then it is easy to check that, according to the
discussion preceding~\eqref{e:triangle}, the homomorphism $F_{-X}$
fits into an exact triangle
\begin{equation}\label{e:triangle1}
\begin{graph}(6,2.3)
\graphlinecolour{1}\grapharrowtype{2}
\textnode {A}(1,1.8){$\hf (S^3)$}
\textnode {B}(5,1.8){$\hf (S^3_{-2g_s} (\mk))$}
\textnode {C}(3,0.3){$\hf (S^3_{-2g_s+1}(\mk))$}
\diredge {A}{B}[\graphlinecolour{0}]
\diredge {B}{C}[\graphlinecolour{0}]
\diredge {C}{A}[\graphlinecolour{0}]
\freetext (2.8,2.1){$F_{-X}$}
\freetext (4.6,0.9){$F_U$}
\freetext (1.4,0.9){$F_V$}
\end{graph}
\end{equation}
where $\overline K$ denotes the mirror image of $K$. By
Theorem~\ref{t:+1} the map $F_{-X}$ sends the non--zero
contact Ozsv\'ath--Szab\'o invariant $c(S^3, \xi _{st})$ to
$c(Y_1(L),\xi_1)$. It is now easy to see that the cobordism $V$ viewed
up--side down is obtained by attaching a 2--handle to $S^3$ along $K$
with framing $2g_s(K)-1$. Therefore, $V$ contains a smoothly
embedded surface of genus $g_s(K)$ and self--intersection
$2g_s(K)-1$. It follows by Proposition~\ref{p:adjunction} that
$F_V=0$. By exactness this means that $F_{-X}$ is injective, therefore
\[
F_{-X}(c(S^3, \xi_{st}))=
c(Y_1(L), \xi_1)\neq 0,
\]
and the claim~\eqref{e:invariant} is proved for $k=1$. We are left to
prove that 
\begin{equation}\label{e:step}
c(Y_{k},\xi_{k})\neq 0\qua\Longrightarrow\qua c(Y_{k+1},\xi_{k+1})\neq 0
\end{equation}
for every $k\geq 1$.

By construction, $(Y_{k+1},\xi_{k+1})$ is given as contact
$(+1)$--surgery on a Legendrian knot in $(Y_{k},\xi_{k})$.  If
$X_k$ denotes the corresponding cobordism, by Theorem~\ref{t:+1} we
have
\begin{equation}\label{e:image}
F_{-X_k}(c(Y_{k}, \xi_{k}))= c(Y_{k+1},\xi_{k+1}). 
\end{equation}

\begin{figure}[ht]
\setlength{\unitlength}{1mm}
\begin{center}
\begin{picture}(125,76)
\put(0,0){\includegraphics[height=7.6cm]{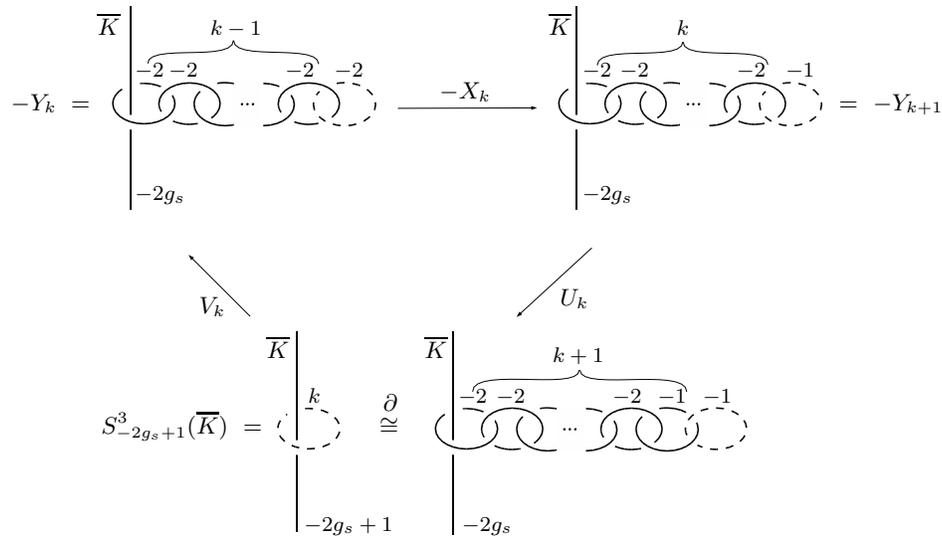}}
\put(3,60){\footnotesize $-Y_{k}\ =\ $}
\put(60,61.5){\footnotesize $-X_k$}
\put(14.5,70){\footnotesize $\overline{K}$}
\put(19.5,64.5){\scriptsize $-2$}
\put(19.5,48){\scriptsize $-2g_s$}
\put(24,64.5){\scriptsize $-2$}
\put(39.5,64.5){\scriptsize $-2$}
\put(46,64.5){\scriptsize $-2$}
\put(29.5,70){\scriptsize $k-1$}
\put(74.5,70){\footnotesize $\overline{K}$}
\put(79,64.5){\scriptsize $-2$}
\put(79,48){\scriptsize $-2g_s$}
\put(84,64.5){\scriptsize $-2$}
\put(99.5,64.5){\scriptsize $-2$}
\put(106,64.5){\scriptsize $-1$}
\put(91.5,70){\scriptsize $k$}
\put(112,60){\footnotesize $\ =\ -Y_{k+1}$}
\put(76,34){\footnotesize $U_k$}
\put(28,33){\footnotesize $V_k$}
\put(37,27){\footnotesize $\overline{K}$}
\put(42,4){\scriptsize $-2g_s+1$}
\put(15,17){\footnotesize $S^3_{-2g_s+1}({\overline K}) \ =\ $}
\put(52,17){\footnotesize $\cong$}
\put(52.5,20){\footnotesize $\del$}
\put(58,27){\footnotesize $\overline{K}$}
\put(62.5,21){\scriptsize $-2$}
\put(63,4){\scriptsize $-2g_s$}
\put(67.5,21){\scriptsize $-2$}
\put(83,21){\scriptsize $-2$}
\put(89,21){\scriptsize $-1$}
\put(95,21){\scriptsize $-1$}
\put(75,26.5){\scriptsize $k+1$}
\put(42.5,21){\scriptsize $k$}%
\end{picture}    
\end{center}
\caption{The surgery exact triangle involving 
$-Y_{k}$, $-Y_{k+1}$ and $S^3_{-2g_s+1}({\overline K})$.}
\label{f:trian}
\end{figure}

The homomorphism $F_{-X_k}$ fits into the
exact triangle
\begin{equation}\label{e:triangle3}
\begin{graph}(6,2)
\graphlinecolour{1}\grapharrowtype{2}
\textnode {A}(1,1.5){$\hf (-Y_{k})$}
\textnode {B}(5, 1.5){$\hf (-Y_{k+1})$}
\textnode {C}(3, 0){$\hf (S^3_{-2g_s+1} ({\overline {K}}))$}
\diredge {A}{B}[\graphlinecolour{0}]
\diredge {B}{C}[\graphlinecolour{0}]
\diredge {C}{A}[\graphlinecolour{0}]
\freetext (3,1.8){$F_{-X_k}$}
\freetext (4.6,0.6){$F_{U_k}$}
\freetext (1.4,0.6){$F_{V_k}$}
\end{graph}
\end{equation}
where $\overline K$ denotes the mirror image of $K$ and the cobordisms
$-X_k$, $U_k$ and $V_k$ are described in Figure~\ref{f:trian} where,
in each picture, the framed dashed knot represents the attaching
circle of a 2--handle giving rise to a cobordism. Remarkably, the
third manifold in the triangle is independent of $k$. This is evident
from the diffeomorphism given in the lower portion of
Figure~\ref{f:trian}, which is obtained by $k+1$ blowdowns.
We are going to show that, for every $k\geq 1$, the cobordism $V_k$
contains an embedded surface $\Si$ of genus $g(\Si)>0$ and
\begin{equation}\label{e:adj-violated}
\Si\cdot\Si \geq 2g(\Si) - 1.
\end{equation}
In view of Proposition~\ref{p:adjunction}, this implies $F_{V_k}=0$,
and therefore that $F_{-X_k}$ is injective. Assuming
$c(Y_{k},\xi_{k})\neq 0$, Equation~\eqref{e:image} then implies
$c(Y_{k+1},\xi_{k+1})\neq 0$, and~\eqref{e:step} follows. Therefore, to finish
the proof we only need to establish the existence of the surface
$\Si\subset V_k$ satisfying~\eqref{e:adj-violated}.

The cobordism $V_k$ is obtained by attaching a 2--handle to
$S^3_{-t}({\overline K})$, where the corresponding framed attaching
circle is shown in the lower left portion of Figure~\ref{f:trian}. We
can think of $S^3_{-t}({\overline K})$ as the boundary of the
4--manifold $Z$ obtained by attaching a 2--handle $H_{\overline{K}}$
to the 4--ball along $\overline{K}$ with framing $-t$. Let $W$ denote
the union $Z\cup V_k$, and let $F\subset Z$ be a smooth surface
representing a generator of $H_2(Z;\Z)$ obtained by capping off a
slicing surface for $\overline{K}$ with the core disk of
$H_{\overline{K}}$. Consider a generic pushoff $F'$ of $F$, viewed as
a surface in $W$. When suitably oriented, $F$ and $F'$ intersect
transversely in $t$ negative points $p_1,\ldots, p_{t}\in F'$.
Consider $t$ generic pushoffs $S_1,\ldots, S_t$ of the embedded
2--sphere $S\subset W$ corresponding to the $k$--framed unknot of the
lower left portion of Figure~\ref{f:trian}, oriented so that $S_i\cdot
F=+1$ for $i=1,\ldots, t$. Each 2--sphere $S_i$ intersects $F$
transversely in a unique point $q_i$. Consider disjoint, smootly
embedded arcs $\ga_1,\ldots,\ga_t\subset F$ such that $\ga_i$ joins
$p_i$ to $q_i$ for each $i=1,\ldots,t$. Let $\nu(F)$ be a small
tubular neighborhood of the surface $F$.  We can view its boundary
$\del\nu(F)$ as a smooth $S^1$ bundle
\[
\pi\co\del\nu(F)\to F, 
\]
so that each of the sets $F'\cap\del\nu(F)$ and $\cup_{i=1}^t S_i
\cap\del\nu(F)$ consists of exactly $t$ fibers of $\pi$. The immersed surface
\[
\widetilde\Si = F'\setminus\nu(F)\cup_{i=1}^t\pi^{-1}(\ga_i)
\cup_{i=1}^t S_i\setminus\nu(F)\subset W
\]
is contained in the complement of $F$. The singularities of
$\widetilde\Si$ come from the intersections among $S_1,\ldots,S_t$ and
$F'$. Resolving those singularities one gets a smoothly embedded
surface which can be isotoped to a surface $\Si\subset V_k$.
Moreover, a simple computation using the fact that $g(F')=g_s(K)=\frac
12(t+1)$ shows that
\[
\Si\cdot\Si = t^2k+t\quad\text{and}\quad 
g(\Si) = \frac {t(t-1)}2 k + \frac{t+1}2.
\]
Since 
\[
\Si\cdot\Si - (2g(\Si)-1) = tk>0,
\]
the surface $\Si$ satisfies~\eqref{e:adj-violated}. This concludes the
proof of Theorem~\ref{t:main}.\qed

\section{The proof of Theorem~\ref{t:nonfill}}\label{s:nonfill}

The following facts~\eqref{e:lens},~\eqref{e:reverse-orientation}
and~\eqref{e:inequality} are proved in~\cite[Propositions~3.1 and
5.1]{OSzF2}.  Let $L(p,q)$ be a lens space. Then,
\begin{equation}\label{e:lens}
\dim _{\Z/2\Z}\hf (L(p,q))=p.
\end{equation}
Let $Y$ be a closed, oriented 3--manifold, and let $-Y$ be the same
3--manifold with reversed orientation. Then,
\begin{equation}\label{e:reverse-orientation}
\hf(-Y)\cong \hf (Y).
\end{equation}
If $b_1(Y)=0$ then 
\begin{equation}\label{e:inequality}
\dim _{\Z/2\Z}\hf (Y)\geq \vert H_1(Y; \Z)\vert. 
\end{equation}

A rational homology 3--sphere $Y$ is called an \emph{$L$--space} if
\[
\dim _{\Z/2\Z} \hf (Y) = \vert H_1(Y; \Z )\vert.
\]
Notice that according to~\eqref{e:lens} lens spaces are $L$--spaces.

\begin{prop}\label{p:lsp}
Let $K\subset S^3$ be a knot such that $g_s(K)>0$ and $S^3_n(K)$ is an
$L$--space for some integer $n>0$. Then, $S^3_r(K)$ is an $L$--space
for every rational number $r\geq 2g_s(K)-1$.
\end{prop}

\begin{proof}
The 3--manifold $S^3_r(K)$ is an $L$--space for every rational
number $r\geq n$. In fact, it follows
from~\cite[Proposition~2.1]{OSzlens}, that
\begin{equation}\label{e:ab}
S^3_{\frac ab}(K)\qua\text{$L$--space}\qua
\Longrightarrow\qua S^3_{\frac{a+1}b}(K)\qua
\text{$L$--space}.
\end{equation}
Suppose $r=\frac pq\geq n$, and write $p=qn+k$ with $n, k\geq 0$.
Then, applying~\eqref{e:ab} $k$ times starting from
$S^3_{n=\frac{p-k}q}(K)$ one deduces that $S^3_r(K)$ is an $L$--space.

The statement follows immediately if $n<2g_s(K)-1$. If $n\geq
2g_s(K)-1$, it is enough to show that $S^3_{2g_s(K)-1}(K)$ is an
$L$--space. We do this by backwards induction on $n$. For
$n=2g_s(K)-1$ the statement trivially holds. If $n>2g_s(K)-1$,
consider the surgery exact triangle given by $S^3$ and $K\subset S^3$
with framing $n-1$:
\begin{equation}\label{e:triangle2}
\begin{graph}(6,2)
\graphlinecolour{1}\grapharrowtype{2}
\textnode {A}(1,1.5){$\hf(S^3)\cong\Z/2\Z$}
\textnode {B}(5, 1.5){$\hf (S^3_{n-1} (K))$}
\textnode {C}(3, 0){$\hf (S^3_{n}(K))$}
\diredge {A}{B}[\graphlinecolour{0}]
\diredge {B}{C}[\graphlinecolour{0}]
\diredge {C}{A}[\graphlinecolour{0}]
\freetext (3,1.8){$F_{X}$}
\freetext (4.6,0.6){$F_U$}
\freetext (1.4,0.6){$F_V$}
\end{graph}
\end{equation}
Since the cobordism $X$ contains a smoothly embedded surface $\Sigma$
of genus $g(\Si)=g_s(K)>0$ and 
\[
\Si\cdot\Si=n-1>2g_s(K)-2, 
\]
by Proposition~\ref{p:adjunction} we have $F_X=0$. This implies that
the exact triangle splits, therefore
\[
\hf (S^3_{n}(K))\cong\hf (S^3_{n-1} (K))\oplus\Z/2\Z. 
\]
Hence, if $S^3_{n}(K)$ is an $L$--space then so is $S^3_{n-1} (K)$ once
$n>2g_s(K)-1$, proving the inductive step.
\end{proof}

The following theorem generalizes a result of the first author~\cite{PLpos}:
Recall that  $T_{p,q}$ denotes the positive torus knot of type $(p,q)$. 

\begin{thm}\label{t:nonf}
For each rational number $r\in [2n-1, 4n)\cap\Q$, the 3--manifold 
\[
S^3_r (T_{2n+1,2})
\]
carries no fillable contact structures. 
\end{thm}

\begin{proof}
Figure~\ref{f:plumb} describes a 6--step sequence of 3--dimensional
Kirby moves which show that the oriented 3--manifold
$S^3_r(T_{2n+1,2})$ is the boundary of the 4--dimensional plumbing $X$
described by the last picture.
\begin{figure}[ht]
\begin{center}
\setlength{\unitlength}{1mm}
\begin{picture}(100,100)(0,0)
\put(0,0){\includegraphics[height=10cm]{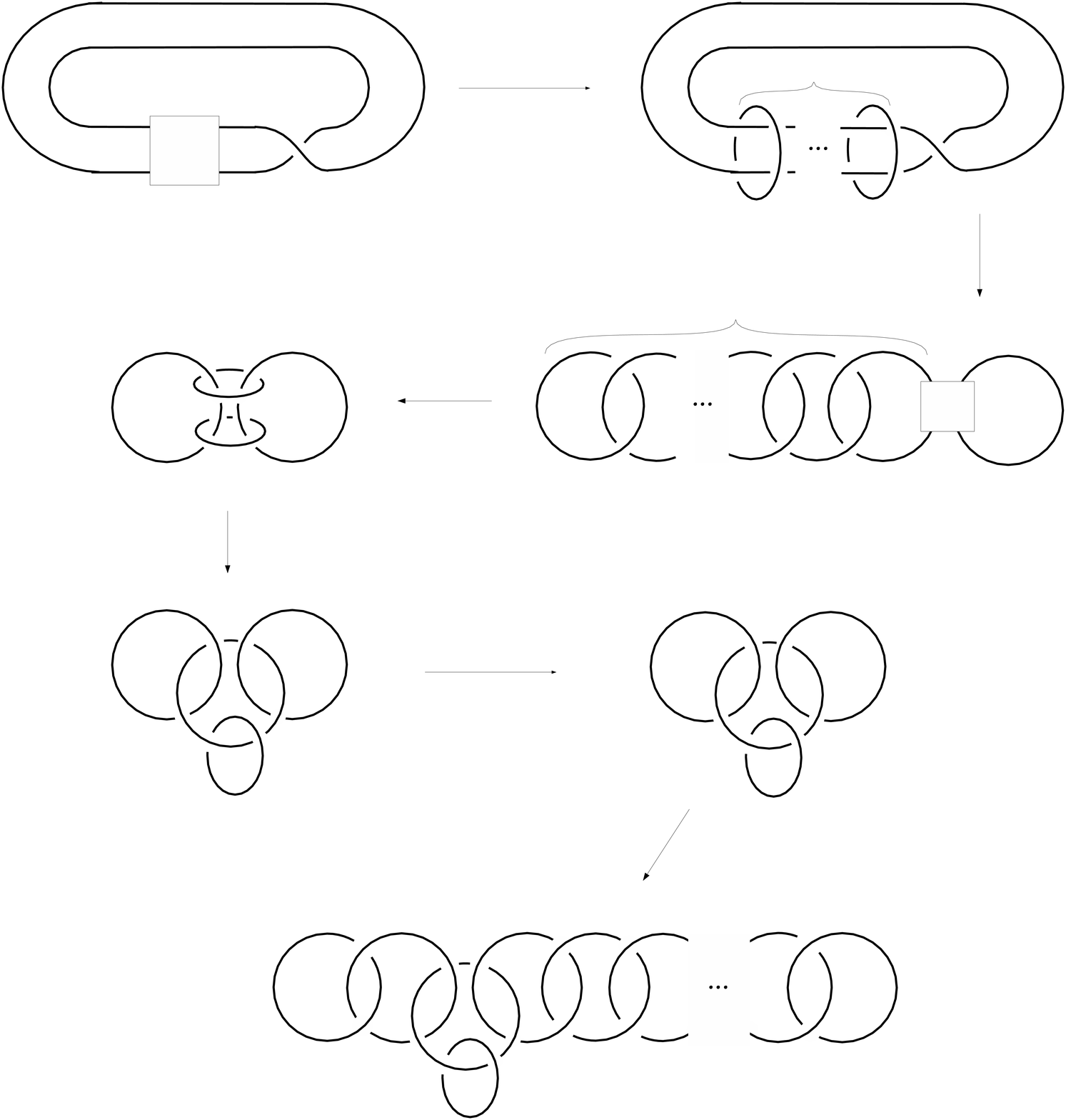}}
\put(15.5,85){\footnotesize $n$}
\put(63,78){\footnotesize $-1$}
\put(73,78){\footnotesize $-1$}
\put(70,92){\footnotesize $n$}
\put(35,97){\footnotesize $r$}
\put(89,97){\footnotesize $r-4n$}
\put(90,68){\footnotesize $r-4n$}
\put(81,62.5){\footnotesize $2$}
\put(62.85,72){\footnotesize $n$}
\put(48.5,56){\footnotesize $-2$}
\put(54,56){\footnotesize $-2$}
\put(62,56){\footnotesize $-2$}
\put(67.5,56){\footnotesize $-2$}
\put(73,56){\footnotesize $-1$}
\put(28,69){\scriptsize $r-4n-2$}
\put(13,65){\scriptsize $-1$}
\put(13,61){\scriptsize $-1$}
\put(5,69){\scriptsize $-\frac{2n+1}n$}
\put(28,47){\scriptsize $r-4n-2$}
\put(12,41){\scriptsize $-1$}
\put(4,46){\scriptsize $-\frac{2n+1}n$}
\put(14,31){\scriptsize $-2$}%
\put(74,47){\scriptsize $\frac{r-4n-2}{r-4n-1}$}%
\put(61,41){\scriptsize $2$}
\put(51,46){\scriptsize $\frac{2n+1}{n+1}$}
\put(63,31){\scriptsize $2$}
\put(25,19.5){\scriptsize $n+1$}
\put(35,19.5){\scriptsize $2$}
\put(45,19.5){\scriptsize $2$}
\put(50,19.5){\scriptsize $a_1$}
\put(56,19.5){\scriptsize $a_2$}
\put(64,19.5){\scriptsize $a_{k-1}$}
\put(74,19.5){\scriptsize $a_k$}
\put(35.5,5.5){\scriptsize $2$}
\put(44,3){\scriptsize $2$}
\put(12,13){$X\ =\ $}
\end{picture}    
\end{center}
\caption{Presentation of $S^3_r(T_{2n+1,2})$ as boundary of a plumbing} 
\label{f:plumb}
\end{figure}
The first step of the sequence is obtained by $n$ blowups. The second
step by $n-1$ handle slides and the third one by two blowups plus a
conversion from integer to rational surgery. The fourth step is given
by a handle slide, the fifth one by three Rolfsen twists and the sixth
one by a conversion from rational to integer surgery. Observe that 
\[
1<\frac{r-4n-2}{r-4n-1}<2
\]
because $r<4n$. The coefficients $a_1,\ldots, a_k$ are given by
\[
\frac{r-4n-2}{r-4n-1} = 
2 - \cfrac{1}{a_1 -
       \cfrac{1}{\ddots -
        \cfrac{1}{a_k}
}},\quad
a_1,\ldots,a_k\geq 2.
\] 
By using~\cite[Theorem 5.2]{NR}, it is easy to check that the 4--dimensional
plumbing $X$ is positive definite. Moreover, the intersection lattice of the
plumbing with reversed orientation $-X$ contains the intersection lattice
$\Lambda_{a_1,n}$ described in Figure~\ref{f:lambda}.
\begin{figure}[ht]
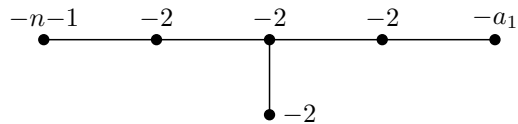

\begin{center}\small
\begin{graph}(6,1.5)(0,0)
\graphnodesize{0.15}
  \roundnode{n1}(0,1)
  \roundnode{n2}(1.5,1)
  \roundnode{n3}(3,1)
  \roundnode{n4}(4.5,1)
  \roundnode{n5}(6,1)
  \roundnode{n6}(3,0)
  
  \edge{n1}{n2}
  \edge{n2}{n3}
  \edge{n3}{n4}
  \edge{n4}{n5}
  \edge{n3}{n6}

  \autonodetext{n1}[n]{$-n{-}1$}
  \autonodetext{n2}[n]{$-2$}
  \autonodetext{n3}[n]{$-2$}
  \autonodetext{n4}[n]{$-2$}
  \autonodetext{n5}[n]{$-a_1$}
  \autonodetext{n6}[e]{$-2$}
\end{graph}
\end{center}
\caption{The intersection lattice $\La_{a_1,n}$}
\label{f:lambda} 
\end{figure}

By~\cite[Theorem~1.4]{OSzUj}, every symplectic filling $(W, \omega)$
of a contact 3--manifold $(Y, \xi)$ such that $Y$ is an $L$--space
satisfies $b_2^+(W)=0$. Since $S^3_{4n+1}(T_{2n+1,2})$ is a lens
space~\cite{Mo} and, by~\cite{KM}, $2g_s(T_{2n+1,2})-1=2n-1$,
Proposition~\ref{p:lsp} implies that $S^3_{r}(T_{2n+1,2})$ is an
$L$--space for every $r\geq 2n-1$. Therefore, every symplectic filling
of a contact 3--manifold of the form $(S^3_r(T_{2n+1,2}),\xi)$ with
$r\geq 2n-1$ satisfies $b_2^+(W)=0$.

If $r\in [2n-1,4n)$, since $Y=S^3_r(T_{2n+1,2})$ is a rational
homology sphere we can build a negative definite closed 4--manifold
\[
Z=W\cup _Y (-X)
\]
which, according to Donaldson's celebrated theorem~\cite{Do1, Do},
must have intersection form $Q_Z$ diagonalizable over $\Z$. Since the
intersection form $Q_{-X}$ embeds in $Q_Z$ it follows that $\La_{a_1,n}$
must embed in $Q_Z$ as well. But we claim that $\La_{a_1,n}$ does not
admit an isometric embedding in the diagonal lattice 
$\Dm=\oplus _m\langle -1\rangle $.
This contradiction forbids the existence of the symplectic filling $W$. 

To prove the claim, we argue as in~\cite[Lemma~3.2]{Paolo}. Suppose
there is an isometric embedding $\varphi$ of $\La_{a_1,n}$ into
$\Dm$. Let $e_1,\ldots, e_k$ be generators of $\Dm$
with self--intersection $-1$. It is easy to check that, up to
composing $\varphi$ with an automorphism of $\Dm$, the four
generators of $\La_{a_1,n}$ corresponding to the vertices of weight
$(-2)$ are sent to $e_1-e_2$, $e_2-e_3$, $e_3-e_4$ and $e_3+e_4$. Up to
composing $\varphi$ with the automorphism of $\Dm$ which sends
$e_4$ to $-e_4$ and fixes the remaining ones, the image $v$ of one of
the two remaining generators of $\La_{a_1,n}$ satisfies
\[
v\cdot (e_3-e_4) = 0,\quad v\cdot (e_3+e_4) = 1,
\]
which is impossible because $(e_3+e_4)-(e_3-e_4)=2e_4$.
\end{proof}

\begin{figure}[ht]
\begin{center}
\setlength{\unitlength}{1mm}
\begin{picture}(65,18)(0,0)
\put(0,0){\includegraphics[height=1.8cm]{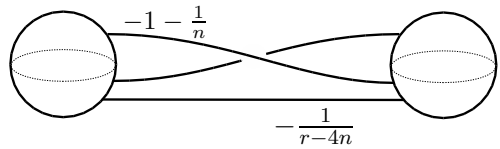}}
\put(15,15){\small$-1-\frac1n$}
\put(35,1){$-\frac1{r-4n}$}
\end{picture}
\end{center}
\caption{Stein fillable contact structures on $Y_{n,r}$ with $r\geq 4n$}
\label{f:stein}
\end{figure}

\rk{Remark} The statement of Theorem~\ref{t:nonf} is optimal, in the sense 
that if $r\not\in [2n-1,4n)$, then the 3--manifold 
\[
Y_{n, r}:=S^3_r (T_{2n+1,2})
\]
supports fillable contact structures. If $r<2n-1$ then, as
observed in the proof of Theorem~\ref{t:main}, $Y_{n,r}$ carries Stein
fillable contact structures. The same holds for $r\geq 4n$. In fact,
examples of Stein fillable contact structures on $Y_{n, r}$ are given by
the contact surgery picture of Figure~\ref{f:stein} (here we are using
our notation as well as the notation of~\cite{G}).

\begin{proof}[Proof of Theorem~\ref{t:nonfill}]
Let $m\in\N$, and let $p_1,\ldots,p_m\in\N$ be consecutive odd primes
with either $p_1=3$ or $p_1=5$, where the choice is made so that
\[
p_1\cdots p_m = 4k+3
\]
for some $k\in\N$. Now let $\al=2k$, and consider the contact
structures obtained via the contact surgeries of
Figure~\ref{f:distinct}.
\begin{figure}[ht]
\begin{center}\small
\setlength{\unitlength}{1mm}
\begin{picture}(45,40)(0,0)
\put(0,0){\includegraphics[height=4cm]{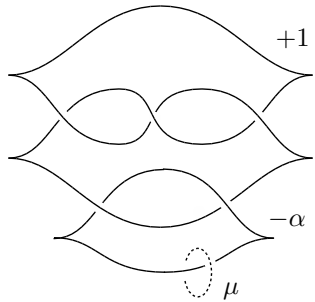}}
\put(37,34){$+1$}
\put(36,10){$-\al$}
\put(30,1){$\mu$}
\end{picture}
\end{center}
\caption{Tight, not fillable contact structures on $N_\al$}
\label{f:distinct}
\end{figure}

The underlying 3--manifold is 
\[
N_\al := S^3_{2+\frac 1{1+\al}}(T_{3,2}).
\]
A simple calculation shows that 
\[
H_1(N_\al;\Z)\cong \Z/(2\al+3)\Z, 
\]
with generator the class of the dotted circle $\mu$ drawn in
Figure~\ref{f:distinct}. The possible choices involved in the contact
surgery construction, ie,~the choices of the $\al-1$ stabilizations
of the Legendrian $(-\al )$--framed unknot, yield contact structures
$\xi_i(\al)$, $i=0,\ldots,\al-1$, where $i$ denotes the number of
right zig--zags added by the stabilizations.  After fixing a suitable
orientation for the knots, this implies that
\[
c_1(\xi_i(\al)) = (2i-(\al-1))\PD([\mu]).
\]
(For computations of homotopic data of contact structures defined by
surgery diagrams see \cite{DGS}.)  Notice that the contact structures
$\xi_i(\al)$ are tight because, since $-\al<0$, they are obtained by
Legendrian surgeries on the contact structure of Figure~\ref{f:knot},
which was shown to have non--zero contact Ozsv\'ath--Szab\'o invariant
in the proof of Theorem~\ref{t:main}. Moreover, since $2+\frac
1{1+\al}\in [1,4)$, by Theorem~\ref{t:nonfill} no $\xi_i(\al)$ is
symplectically fillable.

We claim that, for each $j\in\{1,\ldots,m\}$, there exists an index
$0\leq i(j)<\al$ such that $c_1(\xi_{i(j)}(\al))$ has order
$p_j$. Since the primes $p_j$ are distinct, the claim implies that the
structures $\xi_{i(j)}(\al)$ are pairwise non--isomorphic and, since
$m$ can be chosen arbitrarily large, it suffices to prove the
statement. 

To check the claim, define
\[
i(j) := \frac 12 \left(p_1\cdots\widehat p_j\cdots p_m + \al-1\right).
\]
Then, 
\[
2i(j) - (\al-1) = p_1\cdots\widehat p_j\cdots p_m = \frac{1}{p_j}(2\al +3),
\]
and therefore $c_1(\xi_{i(j)}(\al))$ has order $p_j$. This concludes
the proof.
\end{proof}

\end{document}